\input amstex
\documentstyle{amsppt}
\magnification=\magstep1 \NoRunningHeads

\topmatter
\title
Spectral multiplicities for ergodic flows
\endtitle
\author
Alexandre I. Danilenko and Mariusz Lema{\'n}czyk
\endauthor
\abstract
Let $E$ be a subset of positive integers such that $E\cap\{1,2\}\ne\emptyset$.
 A weakly mixing finite measure preserving flow $T=(T_t)_{t\in\Bbb R}$ is constructed such that the set of spectral multiplicities (of the corresponding Koopman unitary representation generated by $T$) is $E$. Moreover, for each non-zero $t\in\Bbb R$, the set of spectral multiplicities of the transformation $T_t$ is also $E$.
These results are partly extended to actions of some other locally compact second countable Abelian groups.
\endabstract
\thanks
Research of the second named author was
partially supported by Polish MNiSzW grant N N201 384834.
\endthanks
\address
 Institute for Low Temperature Physics
\& Engineering of National Academy of Sciences of Ukraine, 47 Lenin Ave.,
 Kharkov, 61164, UKRAINE
\endaddress
\email alexandre.danilenko\@gmail.com
\endemail

\address
Faculty of Mathematics and Computer Science, Nicolaus Copernicus
University, ul. Chopina 12/18, 87-100 Toru\'n, Poland
\endaddress
\address
 and Institute of Mathematics, Polish Academy of Sciences
ul. \'Sniadeckich 8, 00-950 Warsaw, Poland
\endaddress
\email
mlem\@mat.uni.torun.pl
\endemail

\endtopmatter

\document

\head 0. Introduction
\endhead

Let $G$ be a locally compact second countable Abelian group and let $T=(T_g)_{g\in G}$ be a measure preserving action of $G$ on a standard probability space $(X,\goth B,\mu)$.
The spectral theory of dynamical systems studies the corresponding Koopman unitary representation $U_T=(U_T(g))_{g\in G}$ in the Hilbert space $L_0^2(X,\mu):=L^2(X,\mu)\ominus\Bbb C$ given by
$$
U_T(g)f:=f\circ T_{-g}
$$
(see \cite{KaT}). Such a representation is completely characterized (up to unitary equivalence) by a measure of maximal spectral type on the dual group $\widehat G$ and a spectral multiplicity function $l_T:\widehat G\ni w\to l_T(w)\in\Bbb N\cup\{\infty\}$.
We denote by $\Cal M(T)$ the (essential) image of $l_T$.

One of the most appealing open problems in the spectral theory of dynamical systems can be stated as follows:
when a unitary representation is unitarily equivalent to a Koopman representation?
 Let us consider a  weak version of this problem by replacing the unitary equivalence with another (weaker) equivalence relation on  the set of unitary representations of $G$.
It was introduced in \cite{Fr} for the unitary representations of $\Bbb Z$.
Two unitary representations $U$ and $V$  of $G$ in Hilbert spaces $\Cal H_U$ and $\Cal H_V$ respectively are called {\it cyclicly isomorphic} if there is a unitary operator $W:\Cal H_U\to\Cal H_V$ such that the image under $W$ of each $U$-cyclic subspace in $\Cal H_U$ is a $V$-cyclic subspace in $\Cal H_V$ and vise versa.
Based on the proof in \cite{Fr} for $G=\Bbb Z$ it is easy to see
 that if $U$ and $V$ have a continuous spectrum then they are cyclicly isomorphic if and only if $\Cal M(U)=\Cal M(V)$.
We thus come to the following natural question which is called the spectral multiplicity problem:
\roster
{\it which subsets $E\subset\Bbb N$ are realizable as $E=\Cal M(T)$ for an ergodic (or weakly mixing) free action $T$?}
\endroster

In the case $G=\Bbb Z$ the spectral multiplicity problem was studied by a number of authors (see  references in a recent survey \cite{Le} and subsequent progress in \cite{Ry2}, \cite{KaL},  \cite{Da3}).
It is proved, in particular, that a subset $E\subset\Bbb N$ is realizable if one of the following is fulfilled: $1\in E$, $2\in E$ or $E=n\cdot F$ for some $n\ge 2$ and a subset $F\subset\Bbb N$ with $1\in F$.
It is believed that every subset of $\Bbb N$ is realizable.
In the  case $G=\Bbb Z^2$,  weakly mixing realizations of the subsets $E\ni 1$  were constructed  in \cite{Fi} and weakly mixing realizations of subsets $\{2,4,\dots,2^n\}$, for each $n>0$, were shown  in \cite{Ko}.
For a large class of Abelian locally compact second countable groups $G$ including all countable groups and $\Bbb R^n$, it was proved  in \cite{DaS} that there exist weakly mixing $G$-actions with homogeneous spectrum of arbitrary multiplicity.

In the present paper we mainly consider the case when $G=\Bbb R$.

\proclaim{Theorem 0.1 (Main theorem)}
Let $E$ be a subset of positive integers such that $E\cap\{1,2\}\ne\emptyset$.
\roster
\item"\rom{(i)}"
There is
 a weakly mixing finite measure preserving flow $T=(T_t)_{t\in\Bbb R}$  such that the set of spectral multiplicities of the  Koopman unitary representation generated by $T$ is $E$.
\item"\rom{(ii)}"
For each non-zero $t\in\Bbb R$, the set of spectral multiplicities of the Koopman operator generated by the transformation $T_t$ is also $E$.
\endroster
\endproclaim

Now  Theorem~0.1(i) can be interpreted in the following way: every unitary representation of $\Bbb R$ with continuous spectrum and such that $\Cal M(U)\cap \{1,2\}\ne\emptyset$ is cyclicly isomorphic to a Koopman representation of $\Bbb R$.
Secondly,  given a subset $E\subset\Bbb N$ such that $E\cap\{1,2\}\ne\emptyset$, denote by $\Cal W_E$ the set  of weakly mixing transformations $S$ with $\Cal M(S)=E$.
 As was mentioned above, the set $\Cal W_E$ is known to be non-empty.
 Theorem~0.1(ii) strengthens this fact: $\Cal W_E\cup\{\text{Id}\}$ contains a one-parameter subgroup.

Now we make some remarks concerning the proof of Theorem~0.1.
The simplest way to obtain flows with non-trivial spectral properties is to consider  the suspensions of ergodic transformations with non-trivial spectral properties.
We recall that the suspensions are special flows constructed under the constant function $1$.
In other words, they are $\Bbb R$-actions induced
 by $\Bbb Z$-actions.
 In Section~1 we briefly review  properties of induced actions in a more general setting of pairs $(G,H)$, where $G$ is a locally compact second countable Abelian group and $H$ a closed co-compact subgroup of $G$.
Some of these properties were established in original papers by G.~Mackey \cite{Ma} and R.~Zimmer \cite{Zi}.
By means of the inducing we can obtain ``cheaply''  a realization  of each subset $E\subset\Bbb N$ containing $1$ as the set of  spectral multiplicities of an ergodic flow.
Unfortunately, the condition $1\in E$ is unavoidable within the class of suspension flows.
Moreover, every suspension flow has a non-trivial discrete spectrum.
Therefore to construct weakly mixing realizations we apply another approach.
It is a continuous analogue of the realizations produced in \cite{Da3}.
The desired flows are compact group extensions of either rank-one flows (for realizations of sets $E\ni 1$ in Section~4) or Cartesian squares of rank-one flows (for realizations of sets $E\ni 2$ in Section 5).
Sections~2 and 3 contain some preliminary material to understand the techniques used in Sections 4 and 5.
In the  final Section~6 we partly extend Theorem~0.1 to the actions of other non-compact Abelian groups: connected groups, groups without non-trivial compact subgroups, groups containing a closed one-parameter subgroup, etc.

\head 1. Induced actions
\endhead

Let $G$ be a locally compact second countable Abelian group and $H$ a co-compact subgroup of $G$.
Given a measure preserving action $S=(S_h)_{h\in H}$ of $H$ on a standard probability space $(X,\goth B,\mu)$, we denote by $T=(T_g)_{g\in G}$ the induced action of $G$ on the product space $( G/H\times X ,\lambda_{G/H}\times\mu)$, where $\lambda_{G/H}$ is Haar measure on $G/H$ (see \cite{Ma}, \cite{Zi}).
Fix a Borel cross-section $s:G/H\to G$ of the natural projection
$G\to G/H$ such that $s(H)=0$.
Then
$$
T_g(y,x):=(gy,S_{h(g,y)}x),\tag1-1
$$
where $h(g,y):=-s(gy)+g+s(y)\in H$.
Notice that the  mapping $h:G\times Y\to H$ is a 1-cocycle, i.e.
$$
h(g_1g_2,y)=h(g_1,g_2y)+h(g_2,y)
$$
for all $y\in Y$, $g_1,g_2\in G$.
If $S$ is ergodic then so is $T$.
 Denote by $U_T$ and $U_S$ the Koopman representations of $G$ and $H$ generated by $T$ and $S$ respectively.
Then $U_T$ is unitarily equivalent to the unitary representation of $G$ induced by $U_S$ \cite{Ma}.
Recall that given a unitary representation $V=(V_h)_{h\in H}$
of $H$ in a Hilbert space $\Cal H$, the {\it induced  (by $V$) representation} $U=(U_g)_{g\in G}$ of $G$ is defined on a Hilbert space $L^2(G/H,\lambda_{G/H})\otimes\Cal H$ by the formula
$$
U_{-g}f(y):=V_{h(g,y)}(f(gy)).
$$
Here we consider $f\in L^2(G/H,\lambda_{G/H})\otimes\Cal H$
as a measurable function $f:G/H\to\Cal H$ such that $\int_{G/H}\|f(y)\|^2\,d\lambda_{G/H}(y)<\infty$.
In particular,
under the above identification, if $b\in L^2(G/H,\lambda_{G/H})$
and $a\in{\Cal H}$ then for $f(y)=(b\otimes a)(y)=b(y)a$ and $h\in
H$ we obtain
$$
\left(U_h f\right)(y)=b(y)V_h(a)=(b\otimes V_h a)(y).
$$

\proclaim{Proposition 1.1}
Let $\pi:\widehat G\to\widehat H$ stand for the natural projection.
Denote by  $\sigma_U$  a measure of maximal spectral type of $U$  on $\widehat G$.
\roster
\item"\rom{(i)}"
 $\sigma_U\circ\pi^{-1}$ is a measure of maximal spectral type of $V$.
\item"\rom{(ii)}"
$U$ and $V$ have the same set of spectral multiplicities.
\endroster
\endproclaim
\demo{Proof}
Let $M\subset\Bbb N\cup\{\infty\}$ denote the set of spectral multiplicities of $V$.
Then there is a decomposition $\Cal H=\bigoplus_{i\in M}\bigoplus_{j=1}^i\Cal H_{i,j}$ of $\Cal H$ such that
\roster
\item"$\bullet$" $\Cal H_{i,j}$ is a cyclic space for $V$ for every pair $(i,j)$.
Denote by $\sigma_{i,j}$ a measure of the maximal spectral type for $V\restriction \Cal H_{i,j}$.
Then
\item"$\bullet$"
$\sigma_{i,j}\perp\sigma_{i',j'}$ if $i\ne i'$ and
\item"$\bullet$"
$\sigma_{i,j}\sim\sigma_{i',j'}$ if $i=i'$.
\endroster
It is easy to see that $L^2(G/H,\lambda_{G/H})\otimes\Cal H_{i,j}$ is a cyclic space for $U$.
Denote by $\sigma_{i,j}'$ a measure of the maximal  spectral type of $U\restriction (L^2(G/H,\lambda_{G/H})\otimes\Cal H_{i,j})$.
It is easy to see that
$$
\sigma_{i,j}'\sim\sigma_{i',j'}' \quad\text{if}\quad i=i'.\tag1-2
$$
Let $a\in\Cal H$ be a cyclic vector for $V$ such that the spectral measure of $a$ is $\sigma_V$.
Take a unit vector $b\in L^2(G/H,\lambda_{G/H})$.
Then for each $h\in H$,
$$
\langle U_h(b\otimes a), b\otimes a\rangle=\langle V_ha,a\rangle
$$
This implies that $\pi$ projects the spectral measure of $b\otimes a$ into $\sigma_V$.
This yields $\sigma_{i,j}'\circ\pi^{-1}=\sigma_{i,j}$ for each pair $(i,j)$.
Therefore
$$
\sigma_{i,j}'\perp\sigma_{i',j'}'\quad\text{if}\quad i\ne i'.\tag1-3
$$
Since $L^2(G/H,\lambda_{G/H})\otimes\Cal H=\bigoplus_{i\in M}\bigoplus_{j=1}^iL^2(G/H,\lambda_{G/H})\otimes\Cal H_{i,j}$,
we deduce both (i) and (ii) from \thetag{1-2} plus \thetag{1-3}.
\qed
\enddemo

It follows that if a $G$-action  $T$ is induced by an $H$-action $S$ then $\Cal M(T)=\Cal M(S)\cup\{1\}$.
The ``extra'' value $1$ appears because $L^2_0(G/H,\lambda_{G/H})\otimes 1$ is a $U_T$-cyclic subspace of $L^2_0(G/H\times X,\lambda_{G/H}\times\mu)$.

Denote by $V$ an action of $G$ on the homogeneous  space $G/H$ by translations.
The following  proposition about induced actions was shown by R.~Zimmer in~\cite{Zi}.

\proclaim{Proposition 1.2}
\roster
\item"\rom{(i)}"
Let $T$ be an action of $G$ on $(X,\goth B,\mu)$.
Then the action of $G$ induced by $T\restriction H$ is isomorphic to
the Cartesian product  $V\times T$.
\item"\rom{(ii)}"
An action $T$ of $G$ is induced by an action of $H$ if and only if $T$ has a factor isomorphic to $V$.
\endroster
\endproclaim

Recall that given a dynamical system $(Z,\nu, T)$,  a $T\times T$-invariant measure $\rho$ on the product space $Z\times Z$ such that the  coordinate marginals of $\rho$ are both equal to $\nu$ is called a {\it (2-fold) self-joining of $T$}.
For the theory of joinings and notions like relative weak mixing,
relative compactness, simplicity and  centralizer we refer the
reader to \cite{JuR} and \cite{KaT}.

In the following corollary we describe the structure of self-joinings  of induced actions.

\proclaim{Proposition 1.3}
Let $T$ be a $G$-action induced by an ergodic $H$-action $S$ (see \thetag{1-1} for the notation).
Let $\rho$ be an ergodic  self-joining of $T$.
Then $(Y\times X\times Y\times  X,\rho, T\times T)$ is an induced $G$-action.
More precisely,
there are $\kappa\in J^e_2(S)$ and $g\in G$ such that
$$
\rho=\int_{G/H}\kappa\circ(S_{h(s(y),y)}\times S_{h(s(y),g y)})\times\delta_y\times\delta_{gy}\,d\lambda_{G/H}(y).
$$
 $\rho$ is a graph of an isomorphism if and only if so is $\kappa$.
Hence
two induced $G$-actions are isomorphic if and only if the underlying $H$-actions are isomorphic.
\endproclaim
\demo{Proof}
We use the notation from \thetag{1.1}.
The projection map  $Y\times X\to Y$ intertwines $T$ with $V$ (see \thetag{1.1}).
Therefore the  projection $\rho^*$ of $\rho$ to $Y\times Y$ is an ergodic self-joining of $V$.
Hence there is $g\in G$ such that $\rho^*(A\times B)=\lambda_{G/H}(A\times gB)$ for all measurable subsets $A,B\subset G/H$.
Disintegrate now $\rho$ with respect to $\rho^*$:
$$
\rho=\int_{G/H}\kappa_y\times\delta_y\times\delta_{gy}\,d\lambda_{G/H}(y), \tag1-4
$$
where $Y\times Y\ni(y,y)\mapsto\kappa_y$ is a measurable field of probability measures on $X\times X$ such that
$$
\int_Y\kappa_y^{(1)}\times\delta_y\,d\lambda_{G/H}(y)=
\int_Y\kappa_y^{(2)}\times\delta_{gy}\,d\lambda_{G/H}(y)= \mu\times\lambda_{G/H},\tag1-5
$$
where $\kappa^{(i)}_y$ is the $i$-th coordinate projection of $\kappa_y$ for $i=1,2$ and every $y\in Y$.
Since $\rho$ is $T\times T$-invariant, we deduce from \thetag{1-4} and \thetag{1-1} that
$$
\kappa_{g'y}=\kappa_y\circ(S_{h(g',y)}\times S_{h(g',gy)})\tag1-6
$$
for each $g'\in G$ at a.e. $y\in G/H$.
Substituting  $g'\in H$ into \thetag{1-6} we obtain that
$\kappa_y$ is invariant under $S\times S$ for a.a. $y\in G/H$.
Since $\mu$ is ergodic under $S$, we deduce from~\thetag{1-5} that
$\kappa_y^{(1)}=\kappa_y^{(2)}=\mu$ for a.a. $y\in G/H$.
Thus $\kappa_y$ is a self-joining of $S$ for a.a. $y\in G/H$.
Since $G$ acts transitively on $Y$, the equation \thetag{1-6} can be ``resolved'' in a standard way:
$$
\kappa_{y}=\kappa\circ(S_{h(s(y),y)}\times S_{h(s(y),gy)}), \quad y\in G/H,
$$
for certain self-joining $\kappa$ of $S$ (formally, put $\kappa=\kappa_H$ and $g'=s(y)$ into \thetag{1-6}).
Moreover, $\kappa$ is ergodic.

The
remaining assertions of Proposition 1.3 follow immediately.
\qed
\enddemo

\proclaim{Corollary 1.4}
\roster
\item"\rom{(i)}"
If  $S$ is either relatively weakly mixing or relatively compact with respect to some factor $\goth A$ then $T$ is  either relatively weakly mixing or relatively compact (respectively) with respect to the factor induced by $\goth A$.
\item"\rom{(ii)}"
 $T$ is simple if and only if $S$ has pure point spectrum.
\item"\rom{(iii)}"
$C(T)=\{(\text{\rom{Id}}\times R)T_g\mid g\in G, R\in C(S)\}$.
\item"\rom{(iv)}"
If $\goth F$ is a factor of $T$ that contains the standard factor $V$ then $\goth F$ is an induced action of a factor of $S$.
\endroster
\endproclaim

We note that $T$ may also have  factors which do not contain $V$ (for instance in the case considered in Proposition~1.2(i)).

\comment
Let $\alpha$ be a cocycle of $S$ with values in $K$.
Then the map $\beta:G\times G/H\times X\to K$ given by
$$
\beta(g, y,x):=\alpha(h(g,y),x)
$$
is a cocycle of $T$ with values in $K$.
Of course, the action induced by the $\alpha$-extension of  $S$ is the $\beta$-extension of $T$.
\endcomment

\head 2. Preliminaries
\endhead

We start with an important algebraic lemma.
 Let $G$ be a countable Abelian group, $H$ a subgroup of $G$ and $v:G\to G$ a group automorphism. We set
$$
L(G,H,v):=\{\#(\{v^i(h)\mid i\in\Bbb Z\}\cap H), h\in H\setminus\{0\}\}.
$$

\proclaim{Algebraic Lemma 2.1 (\cite{KwL}, \cite{Da3})}
Given any subset $E\subset\Bbb N$, there exist
a countable Abelian group $G$, a subgroup $H\subset G$ and an automorphism $v:G\to G$ such that
\roster
\item"\rom{(i)}"
$E=L(G,H,v)$ and
\item"\rom{(ii)}"
the subgroup of $\widehat v$-periodic points in $\widehat G$ is countable and dense.
\endroster
\endproclaim

We now recall the definition of rank one.
Let $S=(S_g)_{g\in \Bbb R^d}$ be a measure preserving action of $\Bbb R^d$ on a standard
$\sigma$-finite measure space $(Y,\goth C,\nu)$.

\definition{Definition 2.2}
\roster
\item"(i)"
A {\it Rokhlin tower or column} for $S$ is a triple $(A,f, F)$,
 where $A\in\goth C$, $F$ is a cube in $\Bbb R^d$
and $f:A\to F$ is a measurable mapping
 such that for any Borel subset $H\subset F$ and an element $g\in \Bbb R^d$ with
 $g+H\subset F$, one has
 $f^{-1}(g+H)=S_gf^{-1}(H)$.
\item"(ii)"
 $S$ is said to be of {\it rank one (by cubes)} if
there exists a sequence of Rokhlin towers $(A_n,f_n,F_n)$ such that
the volume of $F_n$ goes to infinity and
for any subset $C\in\goth C$ of finite measure, there is a sequence
of Borel subsets $H_n\subset F_n$ such that
$$
\lim_{n\to\infty}\nu(C\triangle f_n^{-1}(H_n))=0.
$$
\endroster
\enddefinition

\comment
\proclaim{Lemma}
If $U(a_n)\to 0.5(I+ U(-1))$ then $U$ is weakly mixing.
\endproclaim

\proclaim{Lemma}
Let $\Delta_n$ be a sequence of subsets in $X$ and let $A_n$ be a sequence of finite subsets in $\Bbb R$ such that $T_a\Delta_n$, $a\in A_n$, are disjoint.
If for each vector $h$ in $L^2(X,\mu)$ we have
$$
\text{\rom{dist}}(h,\text{\rom{Lin}}\{1_{T_a\Delta_n}\mid a\in A_n\})\to 0
$$
as $n\to\infty$ than $T$ has a simple spectrum.
\endproclaim
\endcomment

Denote by $\Cal R\subset X\times X$ the $T$-orbit equivalence relation. A Borel map $\alpha$ from $\Cal R$ to a compact group $K$ is called a {\it cocycle} of $\Cal R$ if
$$
\alpha(x,y)+\alpha(y,z)=\alpha(x,z)\quad\text{for all }(x,y), (y,z)\in\Cal R.
$$
 Two cocycles $\alpha,\beta:\Cal R\to K$ are {\it cohomologous} if there is a $\mu$-conull subset $B\subset X$ such that
$$
\alpha(x,y)=\phi(x)+\beta(x,y)-\phi(y)\quad
\quad\text{for all }(x,y) \in\Cal R\cap  (B\times B).
$$
for a Borel map $\phi:X\to K$.
Given a cocycle $\alpha:\Cal R\to K$ and a closed subgroup $H\subset K$,
we can define a new flow $T^{\alpha,H}=(T^{\alpha,H}_t)_{t\in\Bbb R}$ on the space $(X\times K/H,\mu\times\lambda_{K/H})$ by setting
$$
T^{\alpha,H}_t(x,k+H)=(T_tx,\alpha(T_tx,x)+k+H).
$$
This flow is called {\it a compact group extension} of $T$.
Given a character $\chi\in\widehat K$, we denote by  $U_{T^\alpha,\chi}$
the following unitary representation of $\Bbb R$ in $L^2(X,\mu)$:
$$
(U_{T^\alpha,\chi}(t)f)(x):=\chi(\alpha(T_{-t}x,x))f(T_{-t}x).
$$
There is a natural decomposition of $U_{T^{\alpha,H}}$ into an orthogonal sum
$$
U_{T^{\alpha,H}}=\bigoplus_{\chi\in\widehat{K/H}}U_{T^\alpha,\chi},
$$
where $\widehat{K/H}$ is considered as a subgroup of $\widehat{K}$.

If a transformation $S$ commutes with $T$ (i.e. $S\in C(T)$) then a cocycle $\alpha\circ S:\Cal R\to K$ is well defined by $\alpha\circ S(x,y):=\alpha(Sx,Sy)$. The important cohomology equation on $\alpha$ mentioned in Section 0 can now be stated as follows
$$
\alpha\circ S\text{ \ is cohomologous to \ }  v\circ\alpha\tag 2-1
$$
for some $S\in C(T)$ and a group automorphism $v:K\to K$.

\head 3. $(C,F)$-flows and $(C,F)$-cocycles
\endhead

To prove Main Theorem we will use  the $(C,F)$-construction (see \cite{Da1} and references therein).  We now briefly outline its formalism.
Let two sequences $(C_n)_{n>0}$ and $(F_n)_{n\ge 0}$ of  subsets in $\Bbb R$  be given such that:
\roster
\item"---"
$F_n=[0,h_n)$, $h_0=1$,
\item"---"  $C_n$ is finite, $\# C_n>1$,  $\min C_n=0$,
\item"---" $F_n+C_{n+1}\subset [0,h_{n+1}-1)$,
\item"---" $(F_{n}+c)\cap (F_n+c')=\emptyset$ if $c\ne c'$, $c,c'\in C_{n+1}$,
\item"---" $\lim_{n\to\infty}\frac{h_n}{\#C_1\cdots\# C_n}<\infty$.
\endroster
Let $X_n:=F_n\times C_{n+1}\times C_{n+2}\times\cdots$. Endow this set with the  standard product Borel structure. The  following map
$$
(f_n,c_{n+1},c_{n+2})\mapsto(f_n+c_{n+1},c_{n+2},\dots)
$$
is a Borel embedding of $X_n$ into $X_{n+1}$.
We now set $X:=\bigcup_{n\ge 0} X_n$ and endow it with the  inductive limit standard Borel structure.
Given a Borel subset $A\subset F_n$, we denote by $[A]_n$ the following cylinder: $\{x=(f,c_{n+1},\dots,)\in X_n\mid f\in A\}$.
The family of all cylinders generates the entire $\sigma$-algebra $\goth B$ on $X$.

Let  $\Cal R$ stand for the {\it tail} equivalence relation on $X$: two points $x,x'\in X$ are $\Cal R$-equivalent if there is $n>0$ such that $x=(f_n,c_{n+1},\dots),\ x'=(f_n',c_{n+1}',\dots)\in X_n$ and $c_m=c_m'$ for all $m>n$.
Of course, $\Cal R$ is a Borel subset of $X\times X$.
It is easy to see that there is only one probability (non-atomic) Borel measure $\mu$ on $X$ which is invariant  under $\Cal R$.
This means that every Borel isomorphism of $X$ whose graph is a subset of $\Cal R$ preserves $\mu$.
We note that the restriction of $\mu$ on $X_n$ is an infinite product
$\nu_n\times\kappa_{n+1}\times\kappa_{n+2}\times\cdots$, where $\kappa_i$ is the equidistribution on $C_i$ and $\nu_{n+1}$ is a measure proportional to $\lambda_\Bbb R\restriction F_{n}$.
Hence for each $n\ge 0$ and a subset $A\subset F_n$,
$$
\mu([A]_n)/\mu(X_n)=\lambda_\Bbb R(A)/h_n.
$$
We now isolate  a subset $\widetilde X\subset X$ such that
$$
\widetilde X\cap X_n:=\{x=(f_n,c_{n+1}, c_{n+2},\dots)\in X_n\mid c_k\ne 0\text{ infinitely often}\}.
$$
Then $X_n$ is Borel, $\Cal R$-saturated and $\mu(\widetilde X)=1$.
Now we define a Borel flow  $T=(T_t)_{t\in\Bbb R}$ on $\widetilde X$ by setting
$$
T_t(f_n,c_{n+1},\dots):=(t+f_n,c_{n+1},\dots )\text{ whenever }t+f_n<h_n,\ n>0.
$$
This formula defines $T_t$ partly on $\widetilde X$.
When $n\to\infty$, $T_t$ extends to the entire  $\widetilde X$.
 It is easy to see that the mapping $\widetilde X\times\Bbb R\ni(x,t)\mapsto T_tx\in\widetilde X$ is Borel and $T_{t_1}T_{t_2}=T_{t_1+t_2}$ for all $t_1,t_2\in\Bbb R$.
Moreover, the $T$-orbit equivalence relation  coincides  with $\Cal R\restriction\widetilde X$.
  It follows that $T$ is  $\mu$-preserving.
In what follows we do not distinguish objects (sets, transformations, etc.) if they agree a.e.
That is why  we consider that $T$ is  defined on the entire $X$.

\definition{Definition 3.1} We call $T$ {\it the $(C,F)$-flow} associated with $(C_{n+1},F_n)_{n\ge 0}$.
\enddefinition

It is easy to see that $T$ is of rank one.
Hence it is free and ergodic.

We recall a concept of $(C,F)$-cocycle (see \cite{Da2}). From now on, the group $K$ is assumed Abelian. Given a sequence of maps $\alpha_n:C_n\to K$, $n=1,2,\dots$, we first define a Borel cocycle $\alpha:\Cal R\cap (X_0\times X_0)\to K$ by setting
$$
\alpha(x,x'):=\sum_{n>0}(\alpha_n(c_n)-\alpha_n(c_n')),
$$
whenever $x=(0,c_1,c_2,\dots)\in X_0$,  $x'=(0,c_1',c_2',\dots)\in X_0$ and $(x,x')\in\Cal R$. To extend $\alpha$ to the entire $\Cal R$, we first define a map $\pi:X\to X_0$ as follows. Given $x\in X$, let $n$ be the least positive integer such that $x\in X_n$. Then $x=(f_n,c_{n+1},\dots)\in X_n$. We  set
$$
\pi(x):=(\underbrace{0,\dots,0}_{n+1\text{ times}}, c_{n+1}, c_{n+2},\dots)\in X_0.
$$
Of course, $(x,\pi(x))\in\Cal R$. Now for each pair $(x,y)\in\Cal R$, we let
$$
\alpha(x,y):=\alpha(\pi(x),\pi(y)).
$$
It is easy to verify that $\alpha$ is a well defined cocycle of $\Cal R$ with values in $K$.

\definition{Definition 3.2} We call $\alpha$ {\it the $(C,F)$-cocycle associated with} $(\alpha_n)_{n=1}^\infty$.
\enddefinition

Suppose we have an invertible measure preserving transformation $S$ of $(X,\mu)$ such that $S$ maps bijectively $\Cal R(x)$ on $\Cal R(Sx)$ for $\mu$-a.a. $x\in X$. (This condition holds, for instance, if $S\in C(T)$.)
Then for each cocycle $\alpha:\Cal R\to K$, we can define a cocycle $\alpha\circ S$ by setting
$$
\alpha\circ S(x,y):=\alpha(Sx,Sy), \quad (x,y)\in\Cal R.
$$
Adapting the argument from \cite{Da2, Section 4} we obtain the following lemma.

\proclaim{Lemma 3.3} Let $\bar z=(z_n)_{n+1}^\infty$ be a sequence of positive reals. Suppose that
$$
\sum_{n>0}\# (C_n\triangle (C_n-z_n))/\#C_{n}<\infty.
 $$
For each $m>0$, we set
 $$
X^{\bar z}_m:=[0, h_m-z_1-\cdots-z_m)\times\prod_{n>m}(C_{n}\cap(C_n-z_n))\subset X_m.
$$
Then a transformation $S_{\bar z}$ of $(X,\mu)$ is well defined by setting
$$
S_{\bar z}(x):=(z_1+\cdots+z_m+f_m, z_{m+1}+c_{m+1}, z_{m+2}+c_{m+2},\dots)\tag3-1
$$
for all $x=(f_m,c_{m+1},c_{m+2},\dots)\in X^{\bar z}_m$, $m=1,2,\dots$.
Moreover, $S_{\bar z}$ commutes with $T$ and $T^{z_1+\cdots+z_m}\to S_{\bar z}$ weakly as $m\to\infty$.

Let $v$ be a continuous group automorphism of $K$ and let
$$
C_m^\circ:=\{c\in C_m\cap (C_m-z_m)\mid \alpha_m(c+z_m)=v(\alpha_m(c))\}. $$
If
$$
\sum_{n>0}(1-\# C_n^\circ/\#C_{n})<\infty\tag3-2
$$
then the cocycle $\alpha\circ S_{\bar z}$ is cohomologous to $v\circ\alpha$.
\endproclaim

\head 4.
 Realization of sets containing $1$ as spectral multiplicities
\endhead

Let $E$ be a subset of positive integers.
 By Algebraic Lemma~2.1, there exist a compact Polish Abelian group $K$,  a closed subgroup $H$ of $K$ and  a continuous automorphism $v$ of $K$ such that
$$
E=L(\widehat K,\widehat{K/H}, \widehat v).
 $$
The subgroup  of $v$-periodic points in $K$ will be denoted by $\Cal K$.
It is countable and dense in $K$ by Lemma~2.1.
Let $\xi_1$ and $\xi_2$ be two rationally independent positive reals in $\Bbb R$.
 Fix a partition
$$
\Bbb N=\Cal W_1\sqcup\Cal W_2\sqcup\bigsqcup_{a\in \Cal K} \Cal N_a
$$
of $\Bbb N$ into infinite subsets.
To construct the desired realization we define inductively a   sequence $(C_n,h_n,\alpha_n)_{n=1}^\infty$, where $C_n$ a finite subset of $\Bbb R$, $h_n$ is a  positive real and $\alpha_n:C_n\to\Bbb R$ a mapping.
Suppose we have already constructed this sequence up to index $n$.
Consider  two cases.

If  $n+1\in \Cal N_a$ for some $a\in\Cal K$ then we denote by $m_a$ the least positive period of $a$ under $v$.
Now we set
$$
\gather
z_{n+1}:=m_anh_n, \quad r_n:=n^3m_a,\\
C_{n+1}:=h_n\cdot\{0,1,\dots, r_n-1\},\\
h_{n+1}:= r_nh_n+1,
\endgather
$$
Let $\alpha_{n+1}:C_{n+1}\to K$ be any map satisfying the following conditions
\roster
\item"(A1)" $\alpha_{n+1}(c+z_{n+1})=v\circ\alpha_{n+1}(c)$ for all $c\in C_{n+1}\cap (C_{n+1}-z_{n+1})$,
\item"(A2)" for each $0\le i< m_a$ there is a subset $C_{n+1,i}\subset C_{n+1}$ such that
$$
\gather
C_{n+1,i}-h_n\subset C_{n+1}, \\  \alpha_{n+1}(c)=\alpha_{n+1}(c-h_n)+v^i(a)\text{ for all }c\in C_{n+1,i}\text{  and}\\
\bigg|\frac{\#C_{n+1,i}}{\#C_{n+1}}-\frac 1{m_a}\bigg|<\frac 2{nm_a}.
\endgather
$$
\endroster

If  $n+1\in \Cal W_i$ for $i=1,2$ then we set
$$
\gather
C_{n+1}:=\{jh_n\mid 0\le j< n\}\sqcup
\{j(h_n+\xi_i)+nh_n\mid 0\le j< n\},
\\
h_{n+1}:= 2nh_n+n\xi_i,\\
\alpha_{n+1}(c):=1_K\quad\text{for all }c\in C_{n+1}.
\endgather
$$

Thus, $C_{n+1},h_{n+1},\alpha_{n+1}$ are completely defined.

 Denote by $(X,\mu,T)$ the  $(C,F)$-flow associated with the sequence $(C_{n+1},F_n)_{n\ge 0}$, where $F_n:=[0,h_n)$. Let $\Cal R$ stand for the tail equivalence relation (or, equivalently, the $T$-orbit equivalence relation) on $X$.
 Denote by $\alpha:\Cal R\to K$ the cocycle of $\Cal R$ associated with the sequence $(\alpha_n)_{n>0}$.
Let $\lambda_{K/H}$ stand for the Haar measure on $K/H$. We denote by $T^{\alpha,H}$ the following flow on the space
$(X\times K/H,\lambda_{K/H})$:
$$
T_t^{\alpha,H}(x,k+H):=(T_tx,\alpha(T_tx,x)+k+H), \quad t\in\Bbb R.
$$

Our purpose in this section is  to prove the following theorem.

\proclaim{Theorem 4.1} $\Cal M(T^{\alpha,H})=E\cup\{1\}$.
\endproclaim

Since
$$
\sum_{n>0}
\frac{\#(C_n\triangle (C_n-z_n))}{\# C_n}=\sum_{n>0}\frac 2{n^2},
$$
it follows from Lemma 3.3 that a transformation $S_{ \bar z}$ of $(X,\mu)$ is well defined by the formula \thetag{3-1} and $S_{\bar z}\in C(T)$.

 It follows from (A1) and (A3) that \thetag{3-2} is satisfied. Hence by Lemma~3.3,
$$
\text{the cocycle $\alpha\circ S_{\bar z}$ is cohomologous to $v\circ\alpha$.}\tag4-1
$$

We need  more notation.
Given $a\in\Cal K$ and $\chi\in \widehat K$,  let
$$
l_\chi(a):= m_a^{-1}\sum_{i=0}^{m_a-1}\chi(v^i(a)),
$$
where $m_a$ stands for the least positive period of $a$ under $v$.

\proclaim{Lemma 4.2} Let    $\chi\in\widehat K$.
Then
\roster
\item"{\text{(i)}}"
$U_{T^\alpha,\chi}(h_n)\to l_\chi(a)\cdot I$ as $\Cal N_a-1\ni n\to\infty$, $a\in \Cal K$.
\item"{\text{(ii)}}"
$U_{T^\alpha,\chi}(h_n)\to 0.5(I+U_{T^\alpha,\chi}(-\xi_j))$
as $\Cal W_j-1\ni n\to\infty$, $j=1,2$.
\endroster
\endproclaim

\demo{Proof}
We will show only (i). The other claim is shown in a similar way.
Let $n\in\Cal N_a$.
 Take any subset $A\subset F_n$.  We note that $[A]_n=[A+C_{n+1}]_{n+1}$.
Therefore it follows from (A4) that
$$
\align
U_{T^\alpha,\chi}(h_n)1_{[A]_n}(x)&=\sum_{i=0}^{m_{a}-1} \chi(\alpha(x,T_{-h_n}x)) 1_{[A+C_{n+1,i}]_{n+1}}(T_{-h_n}x)+ \bar o(x)\\
&=\sum_{i=0}^{m_{a}-1} \chi(v^i(a)) 1_{[A+C_{n+1,i}+h_n]_{n+1}}(x)
+\bar o(x)
\endalign
$$
where $\bar o(x)$ is a function whose $L^2$-norm is small.
Hence
$$
U_{T^\alpha,\chi}(h_n)-\sum_{i=0}^{m_{a}-1} \chi(v^i(a)) 1_{[C_{n+1,i}+h_n]_{n+1}}\to 0
$$
weakly as $\Cal N_{a}-1\ni n\to\infty$, where the  function
$
1_{[C^a_{n+1,i}+h_n]_{n+1}}\in L^\infty(X,\mu)
$
is considered as a multiplication operator in $L^2(X,\mu)$.

It remains to use the inequalities from (A2) and a standard fact that for any sequence $C_{n}'\subset C_{n}$ such that $\# C_{n}'/\# C_n\to\delta$  for some $\delta>0$ we have
$$
1_{[C_{n}']_{n}}\to \delta I \quad\text{weakly as }
n\to\infty.
$$
 \qed
\enddemo

\demo{Proof of  Theorem 4.1}
We first verify that $T^{\alpha}$ is weakly mixing.
Let $U_{T^\alpha}(t)f=\exp{(i\lambda t)}f$ for some $f\in L^2(X\times K)$,
$f\ne 0$ and $\lambda\in \Bbb R$.
It follows from Lemma~4.2(ii) that
$$
U_{T^\alpha}(h_n)\to 0.5(I+U_{T^\alpha}(-\xi_j))
$$
and hence
$$
\exp(i{h_n}\lambda)\to 0.5(1+\exp(-i\lambda\xi_j))
$$
as
$\Cal W_j-1\ni n\to\infty$, $j=1,2$.
Therefore $|1+\exp(-i\lambda\xi_j)|=2$ which implies $\exp(-i\lambda\xi_j)=1$ for $j=1,2$.
Since $\xi_1$ and $\xi_2$ are rationally independent, $\lambda=0$.
It remains to show that $T^\alpha$ is ergodic.
If $\chi\ne 1$ then there is $a\in\Cal K$ with $l_\chi(a)\ne 1$.
If $f\in L^2(X,\mu)$ is invariant under $U_{T^\alpha,\chi}$ then Lemma~4.2(i) yields $f=l_\chi(a)f$.
Hence $f=0$.
If $\chi=1$ then $U_{T^\alpha,\chi}=U_T$.
Since $T$ is ergodic, each $U_{T^\alpha,\chi}$-invariant function is constant.
Thus, we have shown that $U_{T^\alpha}$ is weakly mixing.
Hence $U_{T^{\alpha,H}}$ is also weakly mixing.

To show that  $\Cal M(T^{\alpha,H})=E\cup\{1\}$ we consider
 a natural decomposition
of   $U_{T^{\alpha,H}}$ into
 an orthogonal sum
$$
U_{T^{\alpha,H}}=\bigoplus_{\chi\in \widehat {K/H}}U_{T^\alpha,\chi}.
$$
It is enough
 to prove the following:
\roster
\item"(a)"
$U_{T^\alpha,\chi}$ has a simple spectrum for each $\chi$,
\item"(b)"
$U_{T^\alpha,\chi}$ and $U_{T^\alpha,\xi}$ are unitarily equivalent if $\chi$ and $\xi$ belong to the same $\widehat v$-orbit,
\item"(c)" the measures of maximal spectral type of $U_{T,\chi}$ and $U_{T,\xi}$ are mutually singular if $\chi$ and $\xi$ belong to  different $\widehat v$-orbits.
\endroster

For each $\epsilon>0$ and $n>0$, there are a partition of $F_n$ into intervals $\Delta_0,\dots,\Delta_{M_n}$ and reals $t_1,\dots,t_{M_n}$ such that
$\max_j\text{diam}\,\Delta_j<\epsilon$, $\Delta_j=T_{t_j}\Delta_0$
 and the mapping $[\Delta_j]_n\ni x\mapsto \alpha(T_{-t_j}x,x)\in K$ is constant for  each $1\le j\le M_n$.
This implies~(a).

It is straightforward that \thetag{4-1} implies (b).

If $\chi$ and $\eta$ are non-equivalent then there is $a\in\Cal K$ such that $l_\chi(a)\ne l_\eta(a)$.
Moreover, $U_{T,\chi}^{h_n}\to l_\chi(a)I$ and $U_{T,\eta}^{h_n}\to l_\eta(a)I$
as $\Cal N_{a}-1\ni n\to\infty$ by Lemma~4.2(i).
Hence the measures of maximal spectral types of $U_{T,\eta}$ and $U_{T,\chi}$ are mutually singular.
 Thus (c) holds.
\qed
\enddemo

Now we are going to show the following claim.

\proclaim{Proposition 4.3} $\Cal M(T_t)=E$ for each $t\ne 0$.
\endproclaim

For that we need an auxiliary statement from \cite{LeP}.
Given  a Borel measure $\sigma$ on $\Bbb R$,
we let
$
A_\sigma:=\{t\in\Bbb R\mid \sigma *\delta_t
\not\perp\sigma\}.
$
\proclaim{Lemma 4.4(\cite{LeP})} Let $\sigma$ be a finite Borel measure on $\Bbb R$.
If there are an analytic function $a$ on $\Bbb R$ and a sequence of continuous characters $\xi_n\in\widehat{\Bbb R}$ such that $\xi_n\to\infty$ in $\widehat{\Bbb R}$ and $\xi_n\to a$ weakly in $L^2(\Bbb R,\sigma)$ then for each $t_0\in A_\sigma$ there exists $c\in\Bbb C$ with $|c|=1$ and
$
a(t+t_0)=ca(t)\text{ for each }t\in\Bbb R.
$
\endproclaim

\demo{Proof of Proposition 4.3}
Denote by $\sigma_T$ a probability measure of maximal spectral type for $T$.
We first show that $A_{\sigma_T}=\{0\}$.
It follows from Lemma~4.2(ii) that
$$
U_{T^{\alpha,H}}(h_n)\to 0.5(I+U_{T^{\alpha,H}}(-\xi_j))
$$
weakly as $\Cal W_j-1\ni n\to\infty$, $j=1,2$.
 We deduce from this and Lemma~4.4 that for each $t_0\in A_{\sigma_T}$ and $j=1,2$,  there exists a complex number $c_j$ such that
$$
1+\exp(2\pi i\xi_j(t+t_0))=c_j(1+\exp(2\pi i\xi_jt))
$$
for all $t\in\Bbb R$.
This yields  $c_j=1$ and $\exp(2\pi i\xi_jt_0)=1$ for $j=1,2$.
Since $\xi_1$ and $\xi_2$ are rationally independent, $t_0=0$.

Thus if $0\ne t\in\Bbb R$ then $\sigma_T *\delta_t\perp\sigma_T$.
Hence the natural projection $\Bbb R\to \Bbb R/t\Bbb Z$ is one-to-one on a subset of full $\sigma_T$-measure.
This implies that $\Cal M(T^{\alpha,H})=\Cal M(T^{\alpha,H}_t)$.
\qed
\enddemo

\head 5.
 Realization of sets containing $2$ as spectral multiplicities
\endhead

Now let $E$ be a subset of $\Bbb N$ such that $1\not\in E$.
In this section we will realize the set $E\cup\{2\}$.
We first prove a couple of auxiliary lemmata.

\proclaim{Lemma 5.1} Let $T$ be a weakly mixing flow with a simple spectrum.
Let $\xi_1,\xi_2$ be two rationally independent reals.
Suppose that the weak closure WC$(U_T)$ of the group $\{U_T(t)\mid t\in\Bbb R\}$ contains the following operators:
$$
0.5(I+U_T(j\xi_1)),\quad
0.5(I+U_T(\xi_2))\quad\text{and}\quad
 0.5(I+U_T(\xi_2-\xi_1)),\tag5-1
$$
$j=1,2$.
Then the product flow $T\times T:=(T_t\times T_t)_{t\in\Bbb R}$
has a homogeneous spectrum of multiplicity 2 in the orthocomplement to the constants.
\endproclaim

\demo{Proof}
Let $h$ be a cyclic vector for $U_T$.
Denote by $\Cal C$ the closure of the span of 3 vectors $h\otimes h,U_T(\xi_1)h\otimes h$ and $U_T(\xi_2)h\otimes h$.
It follows  from \thetag{5-1} that $\Cal C$ is invariant under the following operators:
$$
\align
&U_T(j\xi_1)\otimes I+I\otimes U_T(j\xi_1),\quad j=1,2,\tag5-2\\
&U_T(\xi_2)\otimes I+I\otimes U_T(\xi_2),\tag5-3\\
&U_T(\xi_2-\xi_1)\otimes I+I\otimes U_T(\xi_2-\xi_1).\tag5-4
\endalign
$$
Slightly modifying the argument from \cite{Ag} and \cite{Ry1},
we deduce from \thetag{5-2} and \thetag{5-3} that
 $$
U_T(n\xi_1)h\otimes h\in\Cal C\quad\text{and}\quad U_T(n\xi_2)h\otimes h\in\Cal C
$$
for all $n\in\Bbb Z$.
Applying \thetag{5-4} to $U(\xi_1)h\otimes h$ we obtain that $U_T(2\xi_1-\xi_2)h\otimes h\in\Cal C$.
Applying  \thetag{5-2} with $j=-2$ step by step infinitely many times
and then with $j=2$ infinitely many times we obtain that
$U(2n\xi_1-\xi_2)h\otimes h\in\Cal C$ for each $n\in\Bbb Z$.
Next, applying \thetag{5-3} to $U(2\xi_1-\xi_2)h\otimes h$, we deduce that $U(2\xi_1-2\xi_2)\in\Cal C$.
Then again apply infinitely many times \thetag{5-2} with $j=-2$ and $j=2$ to obtain $U(2n\xi_1-2\xi_2)h\otimes h\in\Cal C$.
And so on.
Finally, we obtain that
$$
U(2n\xi_1+m\xi_2)h\otimes h\in\Cal C\quad\text{for all }n,m\in\Bbb Z.
$$
Hence $U(t)h\otimes h\in\Cal C$ for all $t\in\Bbb R$.
Since $h$ is cyclic for $U$, it follows that $\Cal H\otimes h\subset\Cal C$ and therefore  $\Cal C=\Cal H\otimes\Cal H$.

 Denote by $m$  the spectral multiplicity function for $(T_t\times T_t)_{t\in\Bbb R}$ and  denote by $\sigma$ the measure of maximal spectral type for $T$.
By the above, $m(\lambda)\le 3$ for $\sigma*\sigma$-a.a. $\lambda\in\Bbb R$.

On the other hand, since $T$ is weakly mixing and $\sigma\times\sigma=\int_\Bbb R\sigma_\lambda\,d\sigma*\sigma(\lambda)$
and $\sigma_\lambda$ is invariant under the flip mapping $\Bbb R^2\ni(y,z)\mapsto (z,y)\in\Bbb R^2$, it follows that
$m(\lambda)\in\{2,4,\dots\}\cup\{\infty\}$.
Hence $m(\lambda)=2$ a.e.
\qed
\enddemo

\proclaim{Lemma 5.2} Let $U$ and $V$ be
unitary representations of $\Bbb R$ with simple spectrum. Assume
that there are  sequences $a_n\to\infty$, $b_n\to\infty$
$a_n'\to\infty$ and $b_n'\to\infty$ such that
\roster
\item"\rom{(i)}"
$U(a_n)\to 0.5(I+U(\xi))$ and $V(a_n)\to 0.5(I+V(\xi))$,
\item"\rom{(ii)}"
$U(b_n)\to 0.5(dI+U(\xi))$ and $V(b_n)\to 0.5(eI+V(\xi))$,
\item"\rom{(iii)}"
$U(a_n')\to 0.5(I+U(\eta))$ and $V(a_n')\to 0.5(I+V(\eta))$ and
\item"\rom{(iv)}"
$U(b_n')\to 0.5(d'I+U(\eta))$ and $V(b_n')\to 0.5(e'I+V(\eta))$
\endroster
for some $\xi,\eta,d,e,d',e'\in\Bbb R$.
If $d\ne e$, $d'\ne e'$ and $\xi$ and $\eta$ are rationally independent then $U\otimes V$ has also a simple spectrum.
\endproclaim
\demo{Proof} Let $v_1$ and $v_2$ be cyclic vectors for $U$ and $V$.
Denote by $\Cal C$ the $U\otimes V$-cyclic subspace generated by $v_1\otimes v_2$.
It follows from (i) and (ii) that
$$
\gather
(I+U(\xi))v_1\otimes(I+V(\xi))v_2\in\Cal C.\\
(dI+U(\xi))v_1\otimes(eI+V(\xi))v_2\in\Cal C.
\endgather
$$
Hence $U(\xi) v_1\otimes v_2+v_1\otimes V(\xi) v_2\in\Cal C$ and
$dU(\xi) v_1\otimes v_2+ev_1\otimes V(\xi) v_2\in\Cal C$.
This implies, in particular that $U(\xi) v_1\otimes v_2\in\Cal C$.
In a similar way, $U({-\xi}) v_1\otimes v_2\in\Cal C$.
Thus, $(U({\xi})\otimes I)\Cal C=\Cal C$.
In a similar way, we deduce from (iii) and (iv) that
 $(U({\eta})\otimes I)\Cal C=\Cal C$.
Hence $(U(n\xi+m{\eta})\otimes I)\Cal C=\Cal C$ for all $n,m\in\Bbb Z$.
Since $\eta$ and $\xi$ are rationally independent, $\Cal C$ is invariant under the unitary representation $U\otimes I$.
It follows that $\Cal C=\Cal H_1\otimes\Cal H_2$. \qed
\enddemo

Let $K,H,v,\Cal K,\xi_1,\xi_2$ be as in the previous section.
We will assume that $\xi_2>\xi_1$ and put $\xi_3:=\xi_2-\xi_1$.
Fix a partition
$$
\Bbb N=\bigsqcup_{i=1}^3\bigsqcup_{a\in \Cal K} \Cal M_{a,i}\sqcup\Cal N_a
$$
of $\Bbb N$ into infinite subsets.
As in the previous section, to construct the desired realization we define inductively a   sequence $(C_n,h_n,\alpha_n)_{n=1}^\infty$, where $C_n$ a finite subset of $\Bbb R$, $h_n$ is a  positive real and $\alpha_n:C_n\to\Bbb R$ a mapping.
Suppose we have already constructed this sequence up to index $n$.
Consider  two cases.

{\it Case 1.}
If  $n+1\in \Cal N_{a}$ for some $a\in\Cal K$  then we denote by $m_a$ the least positive period of $a$ under $v$.
Now we set
$$
\gather
z_{n+1}:=m_anh_n, \quad r_n:=n^3m_a,\\
C_{n+1}:=h_n\cdot\{0,1,\dots, r_n-1\},\\
h_{n+1}:= r_nh_n,
\endgather
$$
Let $\alpha_{n+1}:C_{n+1}\to K$ be any map satisfying the following conditions
\roster
\item"(A1)" $\alpha_{n+1}(c+z_{n+1})=v\circ\alpha_{n+1}(c)$ for all $c\in C_{n+1}\cap (C_{n+1}-z_{n+1})$,
\item"(A2)" for each $0\le i< m_a$ there is a subset $C_{n+1,i}\subset C_{n+1}$ such that
$$
\gather
C_{n+1,i}-h_n\subset C_{n+1}, \\  \alpha_{n+1}(c)=\alpha_{n+1}(c-h_n)+v^i(a)\text{ for all }c\in C_{n+1,i}\text{  and}\\
\bigg|\frac{\#C_{n+1,i}}{\#C_{n+1}}-\frac 1{m_a}\bigg|<\frac 2{nm_a}.
\endgather
$$
\endroster

{\it Case 2.}
If  $n+1\in \Cal M_{a,i}$ for some $a\in\Cal K$  and $i=1,2,3$ then we denote by $m_a$ the least positive period of $a$ under $v$.
Now we set
$$
\gather
z_{n+1}:=m_an(2h_n+\xi_i),\\
D_{n+1}^1:=h_n\cdot\{0,1,\dots, m_an-1\},\\
D_{n+1}^2:=\{j(h_n+\xi_i)+m_nnh_n\mid 0\le j< m_an\},\\
C_{n+1}:=\bigsqcup_{j=0}^{n^2-1}(jz_{n+1}+(D_{n+1}^1\sqcup D_{n+1}^2)),\\
h_{n+1}:= m_an^3(2h_n+\xi_i),
\endgather
$$
Let $\alpha_{n+1}:C_{n+1}\to K$ be any map satisfying the following conditions
\roster
\item"(B1)" $\alpha_{n+1}(c+z_{n+1})=v\circ\alpha_{n+1}(c)$ for each $c\in C_{n+1}^1\cap (C_{n+1}^1-z_{n+1})$,
\item"(B2)" for each $0\le l< m_a$ there is a subset $D_{n+1,l}\subset D_{n+1}^1$ such that
$$
\gather
D_{n+1,l}-h_n\subset D_{n+1}^1, \\  \alpha_{n+1}(c)=\alpha_{n+1}(c-h_n)+v^l(a)\text{ for all }c\in D_{n+1,l}\text{  and}\\
\bigg|\frac{\#D_{n+1,l}}{\#D_{n+1}^1}-\frac 1{m_a}\bigg|<\frac 2{nm_a},
\endgather
$$
\item"(B3)" $\alpha_{n+1}(c)=1_K$ for each $c\in D_{n+1}^2$.
\endroster

Thus, $C_{n+1},h_{n+1},\alpha_{n+1}$ are completely defined.

 Denote by $(X,\mu,T)$ the  $(C,F)$-flow associated with the sequence $(C_{n+1},F_n)_{n\ge 0}$, where $F_n:=[0,h_n)$. Let $\Cal R$ stand for the tail equivalence relation (or, equivalently, the $T$-orbit equivalence relation) on $X$.
 Denote by $\alpha:\Cal R\to K$ the cocycle of $\Cal R$ associated with the sequence $(\alpha_n)_{n>0}$.
 We denote by $T^{\alpha,H}$ the following flow on the space
$(X\times K/H,\lambda_{K/H})$:
$$
T_t^{\alpha,H}(x,k+H):=(T_tx,\alpha(T_tx,x)+k+H), \quad t\in\Bbb R.
$$

The following lemma is an analogue of Lemma~4.2.
It can be proved in a similar way by using (A2), (B2) and (B3).
We leave details to the reader.

\proclaim{Lemma 5.3} Let  $a\in \Cal K$.
Then for each $\chi\in\widehat K$ and $j>0$
\roster
\item"\rom{(i)}"
$U_{T^\alpha,\chi}(h_n)\to l_\chi(a)\cdot I$
 as  $\Cal N_a-1\ni n\to\infty$ and
\item"\rom{(ii)}"
$U_{T^\alpha,\chi}(jh_n)\to 0.5(l_\chi(ja)I
+ U_{T^\alpha,\chi}(-j\xi_i))$ as  $\Cal M_{a,i}-1\ni n\to\infty$.
\endroster
\endproclaim

Our purpose in this section is  to prove the following theorem.

\proclaim{Theorem 5.4} The transformation $T\times T^{\alpha,H}$ is weakly mixing and $\Cal M(T\times T^{\alpha,H})=E\cup\{2\}$.
\endproclaim
\demo{Proof}
To show that  $\Cal M(T\times T^{\alpha,H})=E\cup\{2\}$ we consider
 a natural decomposition
of   $U_{T\times T^{\alpha,H}}$ into
 an orthogonal sum
$$
U_{T\times T^{\alpha,H}}=\bigoplus_{\chi\in \widehat {K/H}}(U_T\otimes U_{T,\chi}).
$$
It is enough
 to prove the following:
\roster
\item"(a)"
$U_T\otimes U_T$ has a homogeneous spectrum 2 in the orthocomplement to the constants,
\item"(b)"
$U_T\otimes U_{T^\alpha,\chi}$ has a simple spectrum if $\chi\ne 0$,
\item"(c)"
$U_T\otimes U_{T^\alpha,\chi}$ and $U_T\otimes U_{T^\alpha,\xi}$ are unitarily equivalent if $\chi$ and $\xi$ belong to the same $\widehat v$-orbit,
\item"(d)" the measures of maximal spectral type of $U_T\otimes U_{T^\alpha,\chi}$ and $U_{T}\otimes U_{T^\alpha,\xi}$ are mutually singular if $\chi$ and $\xi$ are  not on the same $\widehat v$-orbit.
\endroster

It follows from Lemma 5.3(ii) that WC$(U_T)$ contains  operators $0.5(I+U_T(-\xi_1))$,
$0.5(I+U_T(-2\xi_1))$, $0.5(I+U_T(-\xi_2))$ and $0.5(I+U_T(\xi_1-\xi_2))$. Therefore we deduce (a)  from Lemma~5.1.

Fix a nontrivial $\chi\in \widehat K$.
The unitary representation  $U_{T^\alpha,\chi}$ has a simple spectrum (see the proof of Theorem~3.1).
Moreover,
$$
\align
U_{T^\alpha,\chi}(h_n)&\to 0.5(I+U_{T^\alpha,\chi}(-\xi_1)),\quad
U_{T}(h_n)\to 0.5(I+U_{T}(-\xi_1))\\
&\text{as }\Cal M_{0,1}-1\ni n\to\infty\quad \text{and}
\\
U_{T^\alpha,\chi}(h_n)&\to 0.5(I+U_{T^\alpha,\chi}(-\xi_2)),
\quad
U_{T}(h_n)\to 0.5(I+U_{T}(-\xi_2))
\\
&\text{as }\Cal M_{0,2}-1\ni n\to\infty.
\endalign
$$
 by Lemma~5.3(ii).
Since $\chi$ is nontrivial, it follows from    Algebraic Lem\-ma~1.1 that
 there is $a\in\Cal K$ with $l_\chi(a)\ne 1$.
Again by Lemma~5.3(ii),
$$
\align
U_{T^\alpha,\chi}(h_n)&\to 0.5(l_\chi(a)I+U_{T^\alpha,\chi}(-\xi_1)),\quad
U_{T}(h_n)\to 0.5(I+U_{T}(-\xi_1))\\
&\text{as }\Cal M_{a,1}-1\ni n\to\infty\quad \text{and}
\\
U_{T^\alpha,\chi}(h_n)&\to 0.5(l_\chi(a)I+ U_{T^\alpha,\chi}(-\xi_2)),
\quad
U_{T}(h_n)\to 0.5(I+U_{T}(-\xi_2))
\\
&\text{as }\Cal M_{a,2}-1\ni n\to\infty.
\endalign
$$
 Therefore Lemma 5.2 implies (b).

As in the proof of Theorem~4.3 we can define a transformation $S_{ \bar z}$ of $(X,\mu)$  by the formula \thetag{3-1}.
Then $S_{\bar z}\in C(T)$.
It follows from (A1), (B1) and the definition of $C_{n+1}$ that
 \thetag{3-2} is satisfied.
Hence by Lemma~3.3,
the cocycle $\alpha\circ S_{\bar z}$ is cohomologous to $v\circ\alpha$.
 Therefore the unitary representations $U_{T,\chi}$ and $U_{T,\xi}$ are unitarily equivalent whenever $\chi$ and $\xi$ lie on the same orbit of $\widehat v$.
This yields~(c).

To prove (d), we first find $a\in\Cal G$ such that $l_\chi(a)\ne l_\xi(a)$ (see  claim (ii) of Algebraic Lemma).
It follows from Lemma~5.3(i) that
$$
 U_T(h_n)\otimes U_{T^\alpha,\chi}({h_n})\to l_\chi(a)I \text{ \ and \  } U_T({h_n})\otimes U_{T^\alpha,\xi}({h_n})\to l_\xi(a)I
$$
as $\Cal N_a-1\ni n\to\infty$.
This implies  (d).

Finally, since $\Cal M(T\times T^{\alpha,H})\not\ni 1$, it follows that $T\times T^{\alpha,H}$ is weakly mixing.
\qed
\enddemo

\remark{Remark \rom{5.5}} It follows from Lemma~5.3(ii) that $U_{T^{\alpha,H}}(h_n)\to 0.5(I+U_{T^{\alpha,H}}(-\xi_i))$ as $\Cal M_{0,i}-1\ni n\to\infty$, $j=1,2,3$.
As in Proposition~4.3 we can deduce from this fact that
  $\Cal M(T^{\alpha,H})=\Cal M(T^{\alpha,H}_t)$ for each $t\ne 0$.
\endremark

\head 6. Spectral multiplicities for ergodic actions of other groups
\endhead

 The main result of the paper extends partly to actions of some other locally compact second countable Abelian groups $G$.
If $G$ is compact then each ergodic action $T$ of $G$ has a pure point spectrum and $\Cal M(T)=\{1\}$.
Therefore from now on we assume that $G$ is non-compact.

\proclaim{Corollary 6.1} Let $G$ be a torsion free discrete countable Abelian group and let $E$ be a subset of $\Bbb N$ such that $E\cap\{1,2\}\ne\emptyset$.
Then there is a weakly mixing free action $S$ of $G$ such that $\Cal M(S)=E$.
\endproclaim
\demo{Proof} In the case  when $G=\Bbb Z$ see \cite{Da3} and references therein.
Consider now the case when $G\ne\Bbb Z$.
Then there is an embedding $\phi$ of $G$ into $\Bbb R$  such that the subgroup $\phi(\Bbb R)$ is dense in $\Bbb R$.
Indeed, it is well known that $G$ embeds into $\Bbb Q^\Bbb N$.
In turn, the later group obviously embeds into $\Bbb R$.
It remains to note that if an infinite subgroup of $\Bbb R$  is not isomorphic to $\Bbb Z$ then it is dense in $\Bbb R$.

By Theorem~0.1, there is  a weakly mixing action $T$ of $\Bbb R$ such that $\Cal M(T)=E$.
Then the composition $T\circ\phi=(T_{\phi(g)})_{g\in G}$ is a weakly mixing action of $G$ with $\Cal M(T\circ\phi)=\Cal M(T)=E$. \qed
\enddemo

The first claim of the following lemma is, in fact, a slight generalization of Theorem~4.1.
If we replace (relax) ``weak mixing'' in its statement with ``ergodic'' then it follows   from Theorem~4.1 via  Proposition~1.1.

\proclaim{Lemma 6.2}
Let $A$ be a compact second countable Abelian group.
Let  $E$ be a subset of $\Bbb N$ with $1\in E$.
\roster
\item"\rom{(i)}"
There is a weakly mixing free action $W$ of $\Bbb R\times A$
such that $\Cal M(W)=E$.
\item"\rom{(ii)}"
For each torsion free discrete countable Abelian group $G$,
there is a weakly mixing free action $W$ of $G\times A$ such that
$\Cal M(W)=E$.
\endroster
\endproclaim
\demo{Proof}
(i) Let the objects $K,H,v$ be defined exactly as in Section~4.
We now set $K':=K\times A$, $H':=H\times\{0\}$ and $v':=v\times\text{Id}$.
It is straightforward that
$$
L(\widehat {K'},\widehat{K'/H'}, \widehat{v'})=
L(\widehat {K},\widehat{K/H}, \widehat{v})=E.\tag6-1
$$
Moreover, the subgroup of $v'$-periodic points is countable and dense in
$K'$.
We now construct the skew product flow $T^{\alpha',H'}$ in the same way as in Section~4 but with $K',H',v'$ instead of $K,H,v$.
We note that $T^{\alpha',H'}$ acts on the space $(Y,\nu):=(X\times K/H\times A,\mu\times\lambda_{K/H}\times\lambda_A)$.
Denote by $W=(W_{t,a})_{(t,a)\in\Bbb R\times A}$ the  action of the product group $\Bbb R\times A$ on $(Y,\nu)$ generated by $T^{\alpha',H'}$ and the action of $A$ by rotations along the third coordinate.
Then $W$ is free.
Since $T^{\alpha',H'}$ is weakly mixing, so is $W$.
Denote by $U_W$ the corresponding Koopman unitary representation of $\Bbb R\times A$ in $L^2_0(Y,\nu)$.
We have a decomposition
$$
L^2(Y,\nu)=\bigoplus_{\chi\in\widehat{K/H},\eta\in\widehat A}\Cal H_{\chi,\eta},
$$
where $\Cal H_{\chi,\eta}:=L^2(X,\mu)\otimes\chi\otimes\eta$.
We know from Section~4 that  $\Cal H_{\chi,\eta}$ is
a $U_{T^{\alpha',H'}}$-cyclic subspace for each pair $\chi,\eta$.
It is also a $U_W$-cyclic subspace.
The unitary operator $U_W(0,a)$ acts on $\Cal H_{\chi,\eta}$ by multiplying on $\eta(a)$.
Hence if  $\sigma_{\chi,\eta}$ is a measure of maximal spectral type of
$U_{T^{\alpha',H'}}\restriction \Cal H_{\chi,\eta}$ then the measure
 $\sigma_{\chi,\eta}\times\delta_\eta$ on $\widehat R\times\widehat A$ is a measure of maximal spectral type of $U_W\restriction\Cal H_{\chi,\eta}$.
As was shown in Section~4, if $(\chi,\eta)$ and $(\chi',\eta')$ belong to different $v'$-orbits then  $\sigma_{\chi,\eta}\perp\sigma_{\chi',\eta'}$.
It follows that
$\sigma_{\chi,\eta}\times\delta_\eta\perp \sigma_{\chi',\eta'}\times\delta_{\eta'}$.
On the other hand, if $(\chi,\eta)$ and $(\chi',\eta')$ belong to the same
$v'$-orbit then $\sigma_{\chi,\eta}\sim\sigma_{\chi',\eta'}$.
Moreover, $\eta=\eta'$ by the definition of $v'$.
Hence
$\sigma_{\chi,\eta}\times\delta_\eta\sim \sigma_{\chi',\eta'}\times\delta_{\eta'}$.
These facts plus \thetag{6-1} imply that $\Cal M(W)=E$.

(ii) Consider two cases.
If $G$ is not $\Bbb Z$ then  (ii)
follows from (i) in the very same way as Corollary~6.1 follows from Theorem~0.1.
If $G$ is $\Bbb Z$ then we need to modify the proof of the main result from \cite{Da3} (only the case when $E\ni 1$) in the very same way as we modified the proof of Theorem~4.1 in (i).
\qed
\enddemo

 Let $T_1$ and $T_2$ be  probability preserving ergodic
actions of  locally compact second countable Abelian groups $G_1$ and $G_2$ respectively.
Let $T_1\otimes T_2$ stand for the product action $(g_1,g_2)\mapsto T_1(g_1)\times T_2(g_2)$ of the product group $G_1\times G_2$.
It is easy to see that
\roster
\item"$(\bullet)$"
if $T_1$  has a simple spectrum then $\Cal M(T_1\otimes T_2)=\Cal M(T_2)\cup\{1\}$.
\endroster
As far as we know, this fact was first used in \cite{Fi} for $\Bbb Z^2$-actions.

\proclaim{Corollary 6.3} Let $E\ni 1$.
 If one of the following conditions is satisfied
\roster
\item"\rom{(i)}"
 $G$ contains a closed one-parameter subgroup,
\item"\rom{(ii)}"
$G=D\times F$, where $D$ is a torsion free discrete countable Abelian group and $F$ is a locally compact second countable Abelian group
then
\endroster
there is a free weakly mixing action $T$ of $G$ such that $\Cal M(T)=E$.
\endproclaim

\demo{Proof}
(i)
It follows from  \cite{HR, Theorem~24.30} that $G$ is topologically isomorphic to a product
$\Bbb R\times G'$, where $G'$ is a locally compact Abelian group.

Suppose first that $G'$ is non-compact.
 We now claim that there is a weakly mixing free $G'$-action with a simple spectrum.
 To prove this claim we need several standard auxiliary facts which we state here without proof.
\roster
\item"$-$"
 Let $\Cal A$ be the set of all $G'$-actions on a standard probability space $(X,\mu)$.
A $G'$-action is considered as a continuous map from $G'$ to the Polish (in the weak topology) group of all transformations of $(X,\mu)$.
 Then $\Cal A$ endowed with the topology of uniform convergence on the compact subsets in $G'$ is Polish.
\item"$-$"
The conjugacy class of every free $G'$-action is dense in $\Cal A$.
\item"$-$"
The subset of all weakly mixing $G'$-actions and the subset of all
$G'$-actions with a simple spectrum are both  $G_\delta$ in $\Cal A$.
\item"$-$"
There is a   weakly mixing free $G'$-action and there is a free $G'$-action with a simple spectrum.
\endroster
The claim follows from them  via  a  {\it generic argument}.
Now we deduce the assertion of the Corollary 6.3 from Theorem~0.1 and $(\bullet)$.

Consider now the second case when $G'$ is compact.
Then the assertion of
the Corollary 6.3 follows from Lemma~6.2(i).

(ii) is proved in a similar way by replacing the references to Theorem~0.1 and Lemma~6.2(i) with  references to Corollary~6.1 and Lemma~6.2(ii) respectively.
\qed
\enddemo

We note that if $G$ is connected then (i) is satisfied.
  If $G$ has no non-trivial compact subgroups then one of the  conditions
of Corollary 6.3 is satisfied.

We claim that Theorem~0.1, the case $2\in E$, holds true  if we replace $\Bbb R$-actions with actions of groups $G$ which are isomorphic to the product of $\Bbb R^d$ with  torsion free discrete Abelian groups.
However to prove this fact one has to pass all the way of  Section~5
by adjusting all the arguments from there to the case of $\Bbb R^d$-actions.
We leave this routine to the reader.

\enddocument